\documentclass[review]{elsarticle}

\usepackage{bm,amssymb,amsmath,mathrsfs}
\usepackage{amsthm}
\usepackage{lineno}
\usepackage[usenames]{color}

\usepackage{dcolumn}

\newcommand{\bq}{\begin{equation}}
\newcommand{\eq}{\end{equation}}
\newcommand{\bc}{\begin{center}}
\newcommand{\ec}{\end{center}}
\newcommand{\bit}{\begin{itemize}}
\newcommand{\eit}{\end{itemize}}
\newcommand{\ben}{\begin{enumerate}}
\newcommand{\een}{\end{enumerate}}

\theoremstyle{plain}
\newtheorem{theorem}{Theorem}
\newtheorem*{theorem*}{Theorem}
\newtheorem{proposition}[theorem]{Proposition}

\newtheorem{remark}[theorem]{Remark}

\begin{document}

\journal{(internal report CC23-13)}

\begin{frontmatter}

\title{Factorial Moments of the Geometric Distribution of Order $k$}

\author[cc]{S.~R.~Mane}
\ead{srmane001@gmail.com}
\address[cc]{Convergent Computing Inc., P.~O.~Box 561, Shoreham, NY 11786, USA}

\begin{abstract}
  We derive a simple expression for the $r^{th}$ factorial moment $\mu_{(r)}$  
  of the geometric distribution of order $k$ with success parameter $p\in(0,1)$ (and $q=1-p$) in terms of its probability mass function $f_k(n)$.
  Specifically, $\mu_{(r)} = r!f_k((r+1)k+r)/((qp^k)^{r+1})$.  
\end{abstract}

\vskip 0.25in


\end{frontmatter}

\newpage
For a sequence of independent Bernoulli trials with probability $p\in(0,1)$ of success,  
let $N_k$ be the waiting time for the first run of $k$ consecutive successes. 
Also define $q=1-p$.
Then $N_k$ is said to have the \emph{geometric distribution of order $k$}.  
For $k=1$ it is the geometric distribution.
See, e.g., the material in \cite{BalakrishnanKoutras}--\cite{BarryLoBello}.
Its probability generating function (pgf) is \cite{Feller} 
\bq
\label{eq:pgf}
\mathscr{P}(s) = \frac{p^ks^k(1-ps)}{1-s+qp^ks^{k+1}} = \frac{p^ks^k}{1-qs(1+pk+\dots+k^{k-1}s^{k-1})} \,.
\eq
Its moments (and factorial moments) can be obtained from the derivatives of the pgf.
Our interest in the factorial moments.
Denote the $r^{th}$ factorial moment by $\mu_{(r)}=\mathbb{E}[N_k(N_k-1)\dots(N_k-r+1)]$, for $r=1,2,\dots$.
However, the derivatives of the pgf are tedious to compute.
In this note, we derive a simple expression for $\mu_{(r)}$ in terms of the
probability mass function (pmf) $f_k(n) = P(N_k=n)$, via an alternate route.

First we establish some notation and background information.
The pmf $f_k(n)$ satisfies the recurrence \cite{BarryLoBello}
\bq
\label{eq:recrel}
f_k(n) = qf_k(n-1) + pqf_k(n-2) + p^2qf_k(n-3) +\dots + p^{k-1}qf_k(n-k) \,.
\eq
The initial conditions are $f_k(n)=0$ for $n\in[1,k-1]$ and $f_k(k)=p^k$.
The auxiliary polynomial is
\bq
\label{eq:auxpoly}
\mathscr{A}(z) = z^k -qz^{k-1} -qpz^{k-2} -\dots -qp^{k-1} \,.
\eq
Here $z=1/s$ (see eq.~\eqref{eq:pgf}).
Feller's expression \cite{Feller} for the auxiliary polynomial was actually $\mathscr{A}(1/s)$, but we follow \cite{DM3} and employ $z$ below.
Feller \cite{Feller} proved the roots of the auxiliary polynomial are distinct.
We denote the roots by $\lambda_j$, where $j=0,\dots,k-1$.
Feller \cite{Feller} proved that there was one (unique) positive real root (which he termed the ``principal root'') which we denote by $\lambda_0$ below.
All the roots have magnitude less than unity \cite{DM3}.
Dilworth and Mane \cite{DM3} derived expressions for the roots in terms of Fuss-Catalan numbers.
The details can be found in \cite{DM3} and will not be repeated here.
The pmf $f_k(n)$ is given by a sum over the roots as follows (eq.~(15) in \cite{DM3})
\bq
\label{eq:fkn_sum_coeffs}
f_k(n) = \begin{cases}
\displaystyle
\frac{p^k}{k+1}\,\sum_{j=0}^{k-1} \lambda_j^{n-k}\, \frac{\lambda_j-p}{\lambda_j - k/(k+1)} 
&\qquad \displaystyle
p \ne \frac{k}{k+1} \,,
\\
\displaystyle
\frac{p^k}{k+1}\,\bigl(\, 2\lambda_0^{n-k} + \lambda_1^{n-k} + \dots + \lambda_{k-1}^{n-k} \bigr)
&\qquad \displaystyle
p = \frac{k}{k+1} \,.
\end{cases}
\eq
\begin{proposition}
  The $r^{th}$ factorial moment $\mu_{(r)}$ of the geometric distribution of order $k$ 
with success parameter $p\in(0,1)$ (and $q=1-p$) is given as follows, for $r\ge1$.
\bq
\label{eq:facmom_pmf}
\mu_{(r)} = r!\frac{f_k((r+1)k+r)}{(qp^k)^{r+1}} \,.
\eq
\end{proposition}

\newpage
\begin{proof}
  We treat the case $p \ne k/(k+1)$ below.
  The proof is essentially a matter of summing a geometric series and using the property that for any root, $\lambda_j^k(1-\lambda_j) = p^kq$ (see eq.~\eqref{eq:auxpoly} or \cite{DM3}).
From the definition of the factorial moment of a discrete-valued random variable taking values $n=0,1,\dots$,
\bq
\label{eq:facmom_r_1}
\begin{split}
\mu_{(r)} &= \sum_{n=r}^\infty n(n-1)\dots(n-r+1) f_k(n) 
\\
&= \sum_{n=r}^\infty n(n-1)\dots(n-r+1) \biggl(\frac{p^k}{k+1}\sum_{j=0}^{k-1}\lambda_j^{n-k}\,\frac{\lambda_j-p}{\lambda_j-k/(k+1)}\biggr)
\\
&= \frac{p^k}{k+1}\sum_{j=0}^{k-1}\lambda_j^{r-k}\,\frac{\lambda_j-p}{\lambda_j-k/(k+1)} \biggl(\sum_{n=r}^\infty n(n-1)\dots(n-r+1) \lambda_j^{n-r}\biggr) 
\\
&= \frac{p^k}{k+1}\sum_{j=0}^{k-1}\frac{\lambda_j-p}{\lambda_j-k/(k+1)} \frac{r!\lambda_j^{r-k}}{(1-\lambda_j)^{r+1}}
\\
&= r!\frac{p^k}{k+1}\sum_{j=0}^{k-1}\frac{\lambda_j-p}{\lambda_j-k/(k+1)} \frac{\lambda_j^{kr+r}}{\lambda_j^{(r+1)k}(1-\lambda_j)^{r+1}}
\\
&= r!\frac{p^k}{k+1}\sum_{j=0}^{k-1}\frac{\lambda_j-p}{\lambda_j-k/(k+1)} \frac{\lambda_j^{kr+r}}{(qp^k)^{r+1}}  \qquad\qquad\qquad \textrm{(because $\lambda_j^k(1-\lambda_j) = p^kq$)}
\\
&= r!\frac{1}{(qp^k)^{r+1}}\frac{p^k}{k+1}\sum_{j=0}^{k-1}\lambda_j^{(kr+k+r)-k}\,\frac{\lambda_j-p}{\lambda_j-k/(k+1)} 
\\
&= r!\frac{f_k((r+1)k+r)}{(qp^k)^{r+1}} \,.
\end{split}
\eq
\end{proof}
\begin{remark}
  Observe that the proof employs only the formal properties of the roots.
  The explicit expressions for the roots are not required.
\end{remark}
\begin{remark}
  The case $p = k/(k+1)$ yields the same expression as eq.~\eqref{eq:facmom_pmf}.
The details are omitted.  
\end{remark}

\newpage
The value of the pmf $f_k((r+1)k+r)$ can be obtained from the recurrence in eq.~\eqref{eq:recrel},
else by summing the series in eq.~\eqref{eq:fkn_sum_coeffs},
else from combinatorial sums in \cite{PhilippouGeorghiouPhilippou}--\cite{Muselli}.
Muselli \cite{Muselli} published a particularly efficient sum.
In terms of our notation, it is
\bq
\label{eq:pmf_Muselli}
f_k(n) = \sum_{m=1}^{\lfloor\frac{n+1}{k+1}\rfloor} (-1)^{m-1}p^{mk}q^{m-1}
\biggl[\binom{n-mk-1}{m-2} +q\binom{n-mk-1}{m-1}\biggr] \,.
\eq
Hence we can state the following expression for the $r^{th}$ factorial moment.
We set $n=kr+k+r$ in eq.~\eqref{eq:pmf_Muselli}, hence $(n+1)/(k+1)=(k+1)(r+1)/(k+1)=r+1$.
Then
\bq
\label{eq:famcom_use_Muselli}
\mu_{(r)} = \frac{r!}{(qp^k)^{r+1}}\sum_{m=1}^{r+1} (-1)^{m-1}p^{mk}q^{m-1}
\biggl[\binom{(r+1-m)k+r-1}{m-2} +q\binom{(r+1-m)k+r-1}{m-1}\biggr] \,.
\eq
The sums in eq.~\eqref{eq:famcom_use_Muselli} contain terms which equal zero.
For $m=1$, the first binomial coefficient in eq.~\eqref{eq:pmf_Muselli} is zero.
It is also possible for the numerator to be less than the denominator in other binomial coefficients, i.e.~the binomial coefficient equals zero.
An expression for the pmf which does not contain vanishing binomial coefficients was derived by the author as follows
(i) $f_k(n)=0$ for $n<k$, (ii) $f_k(k)=p^k$, (iii) $f_k(n)=qp^k$ for $n\in[k+1,2k]$ and for $n>2k$
\bq
\label{eq:pmf_me}
\begin{split}
  f_k(n) &= qp^k -\sum_{m=2}^{\lfloor\frac{n+1}{k+1}\rfloor} (-1)^m p^{mk}q^{m-1}\binom{n-mk-1}{m-2} -\sum_{m=2}^{\lfloor\frac{n}{k+1}\rfloor} (-1)^m p^{mk}q^m\binom{n-mk-1}{m-1} \,.
\end{split}
\eq
Using eq.~\eqref{eq:pmf_me} yields the following expression for the $r^{th}$ factorial moment.
\bq
\label{eq:famcom_use_me}
\begin{split}
  \mu_{(r)} &= \frac{r!}{(qp^k)^{r+1}}\biggl\{ qp^k -\sum_{m=2}^{r+1} (-1)^m p^{mk}q^{m-1}\binom{(r+1-m)k+r-1}{m-2}
  \\
  &\quad\qquad\qquad\qquad\quad
  -\sum_{m=2}^r (-1)^m p^{mk}q^m\binom{(r+1-m)k+r-1}{m-1}\biggr\} \,.
\end{split}
\eq
\begin{remark}
The mean and variance are given as follows, which involve $f_k(2k+1)$ and $f_k(3k+2)$ respectively.
\begin{align}
  \mu = \mu_{(1)} &= \frac{f_k(2k+1)}{(qp^k)^2}
\nonumber\\
  &= \frac{qp^k(1-p^k)}{(qp^k)^2} = \frac{1-p^k}{qp^k} \,,
\\
\sigma^2 = \mu_{(2)} -\mu_{(1)}^2 +\mu_{(1)}
&= 2!\frac{f_k(3k+2)}{(qp^k)^3} -\mu^2 +\mu
\nonumber\\
&= \frac{2}{q^2p^{2k}} -\frac{2[1+(k+1)q]}{q^2p^k} +\frac{2}{q}
-\frac{(1-p^k)^2}{(qp^k)^2} + \frac{1-p^k}{qp^k} 
\nonumber\\
&= \frac{1}{(qp^k)^2} -\frac{2k+1}{qp^k} -\frac{p}{q^2} \,.
\end{align}
Both expressions agree with those derived by Feller \cite{Feller}.
\end{remark}

\newpage
\subsection*{\bf\color{red}Addendum 12/30/2023}
It was stated above that for $m=1$, the first binomial coefficient in eq.~\eqref{eq:pmf_Muselli} is zero, but this is not completely true.
Consider $n=k$.
A run of $k$ consecutive successes can be obtained in $k$ Bernoulli trials if and only if every trial is a success, hence the probability is $f_k(k)=p^k$.
Setting $n=k$ in eq.~\eqref{eq:pmf_Muselli} yields
\bq
\label{eq:f_neqk}
\begin{split}
f_k(k) &= \sum_{m=1}^1 (-1)^{m-1}p^{mk}q^{m-1} \biggl[\binom{k-mk-1}{m-2} +q\binom{k-mk-1}{m-1}\biggr] 
\\
&= p^k \biggl[\binom{-1}{-1} +q\binom{-1}{0}\biggr] 
\\
&= p^k \,?
\end{split}
\eq
We must interpret $\binom{-1}{-1}=1$ on the grounds that $\binom{i}{i}=1$ for all values of $i$ (including negative $i$) and
$\binom{-1}{0}=0$ on the grounds that $\binom{i}{j}=0$ if $j>i$ (even if $j=0$).
Hence care must be exercised when using eqs.~\eqref{eq:pmf_Muselli} and \eqref{eq:famcom_use_Muselli}.
The binomial coefficients in eqs.~\eqref{eq:pmf_me} and \eqref{eq:famcom_use_me} do not require such subtleties of interpretation.

\end{document}